\documentclass[11pt]{amsart} 

\usepackage{pdfsync}
\usepackage{lmodern}
\usepackage{amsthm,amsmath,amssymb,amscd,graphicx,enumerate, stmaryrd,xspace,verbatim, epic, eepic,color,url}
\usepackage{tikz}
  \usepackage{pgfplots}
\pgfplotsset{compat=1.16}

 \usepackage{tikz-3dplot}
\usepackage{hyperref}
\usepackage{pgf}

   \newtheorem{thm}{Theorem}[subsection]
      \newtheorem*{thm*}{Theorem}

   \newtheorem{prop}[thm] {Proposition}     
   \newtheorem{lemma} [thm]{Lemma}

   \newtheorem{cor} [thm]{Corollary}

   \newtheorem*{conjecture*}{Conjecture}
   \newtheorem*{remark*}{Conjecture}

\theoremstyle{definition}
 
          \newtheorem*{exercise*}{Exercise}   
     
     \newtheorem{example}[thm]{Example}
  
  \newtheorem{defi}[thm] {Definition}

  \newtheorem{remark} [thm]{Remark}

\usepackage{amsmath}
\usepackage{amssymb,amsthm}
\usepackage{amscd}
\usepackage{enumerate}
\input{xy}
\xyoption{all}
  
\newcounter{prg}[section]

\newcommand{\ov}{\overline}

\newcommand{\la}{\longrightarrow}

\newcommand{\A}{\mathbb{A}}
\newcommand{\C}{\mathbb{C}}

\newcommand{\N}{\mathbb{N}}
\newcommand{\PP}{\mathbb{P}}

\newcommand{\R}{\mathbb{R}}

\newcommand{\cE}{\mathcal {E}}

\newcommand{\codim}{\operatorname{codim}}

\newcommand{\mult}{{\operatorname{mult}}}

\begin{document}
 \title{Tropical curves of unibranch points and hypertangency} 
  \author{Lucia Caporaso and Amos Turchet}
  \address[Caporaso]{Dipartimento di Matematica e Fisica\\ Universit\`{a} Roma Tre \\ Largo San Leonardo Murialdo \\I-00146 Roma\\  Italy }\email{lucia.caporaso@uniroma3.it}
  \address[Turchet]{Dipartimento di Matematica e Fisica\\ Universit\`{a} Roma Tre \\ Largo San Leonardo Murialdo \\I-00146 Roma\\  Italy }\email{amos.turchet@uniroma3.it}
  
   \begin{abstract}
   We study   plane curves  meeting at a single unibranch point and show that they must satisfy two equivalent conditions. 
   A geometric one: the tropical curves  associated  to the contact point are  isomorphic.
   An arithmetic  one:  the local invariants   at the contact point   satisfy a simple relation.    
Furthermore, we prove closed formulas for the delta-invariant of a unibranch singularity, and for the dimension of the locus of curves  with an assigned unibranch point. Our work is  motivated by interest in the Lang's exceptional set.
  \end{abstract}


\maketitle
 
\tableofcontents
 
\section{Introduction}
\subsection{Results}

In this paper we study unibranch points of curves, i.e. points which have a single preimage in the normalization of the curve. Our first goal is to study  
under what conditions   two   projective plane curves of degree at least $2$  intersect in exactly one point, and this   point is unibranch for both of them.
This problem  is   also relevant for the study of the  Lang exceptional set appearing in some   conjectures related to arithmetic and hyperbolic geometry, as we shall explain later in  this introduction.

To a unibranch point $q$ of a plane curve $C \subset \PP^2$ of degree  at least $2$
we associate a type made of a pair of integers $(m,n)$, where $m = \mult_q(C)$ is the multiplicity of $C$ at $q$, and $n$ is the multiplicity of the intersection of $C$ with its tangent line $L$ at $q$, i.e. $(L \cdot C)_q = n$; we say that $q$ is a $(m,n)$-point for $C$.  Then,  in Section~\ref{sec:tropical}, we associate a tropical curve, denoted by  $\Gamma^{(m,n)}$, to a unibranch point of type $(m,n)$.
Our   main result  answering the above problem can be stated as follows.

 \begin{thm}
\label{mirrorin}
 Let $B, C\subset \PP^2$ be integral curves 
 of  degree at least $2$ such that  $B\cap C = \{ q\}$; assume that   
 $q$ is an $(m',n')$-point for $B$ and an $(m,n)$-point for $C$.
 Then there is an isomorphism of tropical curves between  $\Gamma^{(m,n)}$ and  $\Gamma^{(m',n')}$;
 equivalently, 
  $ m'/n' = m/n$.
   \end{thm}
See  Theorem~\ref{thm-mn} and Theorem~\ref{mirror}.
 
Tropical curves are metric graphs, i.e. graphs enriched by the assignment of a length (a positive real number) to each edge; see
  \cite{MZ}.
In the classical literature   non-metric graphs,  such as  Enriques diagrams,  have been widely used to study singularities of plane curves; see \cite{EC} or \cite{Wall}. 
Even though such graphs do not quite fit our purpose,   they inspired our construction of the tropical curve  $\Gamma^{(m,n)}$;  see   Subsection~\ref{sec-trop} for more details.

In the second part of Section~\ref{sec-mirror} we investigate the opposite problem, describing the intersection of two curves   at points of the same  (or,  more generally, equivalent)   type.

Finally, our techniques   enable us to obtain   closed formulas for the delta-invariant of a unibranch singularity,  in Proposition~\ref{delta},   and for the dimension of the space of curves   with an assigned unibranch point, in Proposition~\ref{hyptan}.
The proof of these   formulas  occupy the last section of this article.

\subsection{Context}
Our work in this paper is motivated by the algebro-geometric counterparts  of some outstanding conjectures from diophantine and hyperbolic geometry.
To introduce  them, let   $B\subset \PP^2$ be a complex projective curve of degree at least $3$,  having at most nodal singularities.
The so-called algebraic   exceptional set  of $B$, introduced by Lang, is denoted as follows
\begin{equation}
 \label{eq-EB}
\cE(B)=\{C\subset \PP^2:\ C \text{ rational}, \  \  |\nu_C^{-1}(C\cap B)|\leq 2\}  
\end{equation}
 where $\nu_C: C^\nu \to C$ is the normalization, so that $C^\nu\cong \PP^1$ if $C$ is rational. One of the main ideas behind the circle of conjectures of Lang, Vojta, Demailly and Campana, is that $\cE(B)$ should contain all the infinite families of integral points in the quasi-projective surface $\PP^2 \setminus B$, and all the images of nonconstant holomorphic maps $\C \to \PP^2 \setminus B$. 

 Conjecturally, $\cE(B)$ is finite as soon as $\deg B\geq 4$, which corresponds to the pairs $(\PP^2, B)$ being of log general type. More generally, Lang conjectured that the fact that $\cE(B)$ is a proper closed subset of $\PP^2$ should be equivalent to the fact that $(\PP^2,B)$ is of log general type (and this should hold more in general for projective varieties $X$ with a normal crossing divisor $B$). For $\PP^2$  this conjecture is known only for $B$   reducible   with at least three components; the finiteness follows from \cite{CZGm,Wang_etal},  and we refer to \cite{CT} for further references and explicit bounds. It remains   open when $B$ is irreducible or has only two irreducible components, although \cite{Chen,PR,CRY,ATY} show that it holds for a {\underline{very}} general $B$.

 If $C\in \cE(B)$ and $C\cap B=\{q,q'\}$, then  it is clear that $q$ and $q'$ must be    unibranch  points of $C$. Therefore one is naturally led to investigate unibranch points, and curves having high order of contact at such points. In particular when $B$ has (at least) three components,   any $C \in \cE(B)$ will intersect at least one of the components of $B$ in a single unibranch point. This is precisely the setting of Theorem \ref{mirrorin}.

As an example of application of this type of results we mention  that, in our previous paper \cite{CT}, we used as a crucial ingredient a  special case of Theorem \ref{mirrorin} (Theorem 2.2.1 in loc.cit. which assumes $B$ is  smooth at $q$) to show that $\cE(B)$ is empty for a general $B$ with three irreducible components of total degree at least 5.

\subsection{Notation} We work over $\C$.
Throughout the paper, $S$ will be a smooth  projective surface endowed with a birational morphism to $\PP^2$, and
$C\subset S$ 
 an integral projective curve.
 
Given a point   $p\in C$,
we write $\mult_p(C)$ for the multiplicity of $C$ at $p$; if $m=\mult_p(C)$ we say that $p$ is an $m$-fold point .
We denote by $\nu_C:C^{\nu}\to C$ the normalization. If $|\nu_C^{-1}(p)|=1$ we say that $p$ is {\it unibranch}.

   If $C,B\subset S$ are   reduced  curves   with no components in common, and $p\in C\cap B$, we write
$ 
(C\cdot B)_p 
$ 
for their multiplicity of intersection at $p$. 

If $(C\cdot B)_p=\mult_p(C)\mult_p(B)$ we say that $B$ and $C$ are {\it transverse} at $p$. 

If  $|\nu_C^{-1}(B\cap C)|=1$, we say that   $C$ is {\it hypertangent} to    $B$. 
  
  \subsection*{Acknowledgements.} 
  We are grateful to  Pietro Corvaja, Kristin DeVleming, Luca Ferrigno, Nicola Ottolini, Wei Chen and an anonymous referee for useful comments.
  Special thanks to Laura Capuano for showing us Lemma~\ref{contoLaura}.
 LC is partially supported by  PRIN 2017SSNZAW and PRIN 2022L34E7W,  Moduli spaces and birational geometry.
 AT is partially supported by   PRIN  2022HPSNCR: Semiabelian varieties, Galois representations and related Diophantine problems and   PRIN   2020KKWT53: Curves, Ricci flat Varieties and their Interactions, and is a member of the INDAM group GNSAGA. 

 \section{Unibranch points}
 This section is devoted to foundational concepts on unibranch points of curves on surfaces. 

\subsection{Resolving unibranch points}

\begin{defi}
 Let $q\in C\subset S$ be   a unibranch point; set $m=\mult_q(C)$.
 Let $F\subset S$ be  a smooth  integral  curve, $F\neq C$; set $n=(C\cdot F)_q$.
 If $n>m$ we say that $q$ is an {\it  $(m,n)$-point of $C$ with respect to $F$}.
 
 For brevity   we   write:
 ``$q\in C$ is an $(m,n)$-point   w.r.t.  $F$".
 
   
If  $S=\PP^2$ and $F$ is the tangent line to $C$ at $q$, we simply  say that $q$ is an $(m,n)$-point of $C$.
\end{defi}
For example, a $(1,2)$-point of a   curve $C\subset \PP^2$ is a smooth point which is not a flex. 

By definition,  if $q\in C$ is an  $(m,n)$-point  with respect to $F$, then     $F$ and $C$ are not transverse at $q$.

Our surface $S$ is an iterated blow-up of $\PP^2$, hence every point $q\in S$ has an open neighborhood isomorphic to $\A^2$; we will always choose  coordinates in such a neighborhood so that $q$ is the origin. 
We say that a   curve $F\subset S$ such that $q\in F$ is a {\it line locally at } $q$ if there is an open neighborhood 
$U\cong \A^2$ of $q$ such that   the   equation of $F \cap U$ has degree $1$. In particular, if $E\subset S$ is an exceptional divisor, then $E$ is a line locally at any of its points.

For   a unibranch point $q\in C\subset S$
we   introduce  an  iterated sequence of blow-ups.
Set $S_1=S$, let $S_2\to S_1$ be the blow-up at $q$, and $E_2\subset S_2$ the exceptional divisor. We denote by $C^2\subset S_2$ the proper transform of $C$ and by $q^2\in C^2$ the point lying over $q$, i.e. $q^2=C^2\cap E_2$.
 We  use superscripts for strict transforms, and subscripts otherwise, but we denote $C=C^1$ and $q=q^1$ for convenience.
We iterate  to get a 
 sequence of arbitrary length
$$
\ldots C^{h+1}\la C^h\la C^{h-1}\la \ldots  \la C^1=C
$$
where  $C^h\subset S_h$ is the proper transform of  $C^{h-1}$ under the blow-up 
of $S_{h-1}$ at $q^{h-1}\in C^{h-1}$,  and $q^h=C^h\cap E_h$, where $E_h$ is the exceptional divisor of $S_h\to S_{h-1}$.
If    $F\subset S$ is a   curve, we write $F^h\subset S_h$ for its proper transform.

We need the following 
 \begin{lemma}   
\label{lm-k} Let
  $q\in C$  be  an  $(m,n)$-point w.r.t. $L$, where $L$ is a  line locally at $q$; 
  set $k=\lfloor n/m\rfloor$. 
 Then 
 $q^h\in C^h$ is an $(m,n-(h-1)m)$-point w.r.t. $L^h$ for every $h\leq k-1$,  
and $q^k$ is an $m$-fold point.
Moreover
 \begin{enumerate}[(a)]
 \item 
 \label{lm-ka}
if $n/m\neq k$  then $q^k\in C^k$ is an $(m,n-(k-1)m)$-point w.r.t. $L^k$, and
$q^{k+1}\in C^{k+1}$ is an $(n-km, m)$-point w.r.t. $E_{k+1}$;

  \item 
   \label{lm-kb}
 if $n/m=k$  then    $C^k$ is transverse to both $E_k$ and $L^k$.
    \end{enumerate}
\end{lemma}

\begin{proof} We proceed as in the proof of  \cite[Lemma 2.1.1]{CT}, of which the present lemma is a generalization.
We write, as above, $q^1=q\in C^1=C\subset S^1=S$.
We   use    local   coordinates at $q^1\in U\cong \A^2_{x_1,y_1}$  so that $q^1=(0,0)$ and  $L=L^1$ has equation $y_1=0$. As $q^1\in C^1$ is an $(m,n)$-point w.r.t. $L^1$, the  
defining polynomial,  $f_1(x_1,y_1)$, of $C^1$
is 
$$
f_1(x_1,y_1)=a_{0, m}y_1^m+   \sum_{ m+1\leq i+j\leq d} a_{i,j}x_1^iy_1^j 
$$
with 
\begin{equation}
\label{eq-anm}
 a_{0, m}\neq 0, \quad a_{n,0}\neq 0,  \quad   a_{i,0}=0\quad  \forall i<n.    
\end{equation}
We   blow-up at $q^1$  and use local coordinates $(x_2,y_2)$   at $q^2=(0,0)$, with  $x_1=x_2$
and   $y_1=y_2x_2$.  The local equation of  the exceptional divisor $E_2$ is $x_2=0$, and   the local equation of $L^2$ is  $y_2=0$.  Let $f_2(x_2,y_2)=0$ be the local equation of $C^2$, so that 
$$ 
f_2(x_2,y_2)= a_{0,m}y_2^m+\sum_{m+1\leq i+j\leq d} a_{i,j}x_2^{i+j-m}y_2^j.
$$ 
We iterate: for every $h\leq k$ we use use local coordinates $(x_h,y_h)$   at $q^h=(0,0)$, with  $x_{h-1}=x_h$
and   $y_{h-1}=y_hx_h$, so that the local equation of $E_h$ is $x_h=0$ and the local equation of $L^h$ is $y_h=0$.
The defining polynomial  of $C^h$ is
$$
f_h(x_h,y_h)= a_{0,m}y_h^m+\sum_{m+1\leq i+j\leq d} a_{i,j}x_h^{i+j-(h-1)m}y_h^j.
$$
By \eqref{eq-anm}, the  smallest power of $y_h$  appearing  as a summand above is $y_h^m$, 
and the smallest power of 
$x_h$ is $x_h^{n-(h-1)m}$. 
We have
\begin{equation}
 \label{eq-nnm}
 n-(h-1)m=n-hm+m\geq n-km+m\geq n-n+m=m 
\end{equation}
using $h\leq k$ in the first inequality and $k\leq n/m$ in the second.

If  $h\leq k-1$,  the first inequality is strict, hence $n-(h-1)m>m$.
Therefore $f_h(x_h,y_h)$ contains the summand $y_h^m$ but not the summand $x_h^m$. Hence, as $q^h\in C^h$ is unibranch,
all terms  of degree at most $m$ divisible by $x_hy_h$ must vanish. 
Hence $\mult_{q^h}(C^h)=m$ and, as   $y_h=0$ is the local equation of $L^h$, we get
$(L^h\cdot C^h)_{q^h}=n-(h-1)m$,
hence $q^h\in C^h$  is an $(m,n-(h-1)m)$-point w.r.t. $L^h$, as claimed.

Let   $h=k$. If $n/m\neq k$ then $k<n/m$ and the second inequality of \eqref{eq-nnm} is strict. Hence, by the above argument, $q^k\in C^k$ is an $(m,n-(k-1)m)$-point w.r.t. $L^k$.
Now look at $q^{k+1}\in C^{k+1}$, we have
$$
f_{k+1}(x_{k+1},y_{k+1})= a_{0,m}y_{k+1}^m+\sum_{m+1\leq i+j\leq d} a_{i,j}x_{k+1}^{i+j-km}y_{k+1}^j.
$$
As before, $y_{k+1}^m$ is the smallest power of $y_{k+1}$  and $x_{k+1}^{n-km}$ the smallest power of $x_{k+1}$,
but now $n-km<m$. Since $q^{k+1}$ is unibranch, by the same  argument, $\mult_{q^{k+1}}(C^{k+1})=n-km$ and, as $x_{k+1}=0$ is the local equation of $E_{k+1}$, we conclude that $q^{k+1}\in C^{k+1}$ is an $(n-km,m)$-point w.r.t. $E_{k+1}$.
  
 Assume $n/m=k$. Then in \eqref{eq-nnm} we have equality.
  Therefore we have both $x_k^{n-km}=x_k^m$  and   $y_k^m$ appearing  in $f_k$.
Since $q^k$ is unibranch for $C^k$, this is possible only if the homogeneous part, $(f_k)_m$, of degree $m$ of $f_k$,  has form
$ 
(f_k)_m=(\alpha y_k+\beta x_k)^m
$ 
with  $\alpha,\beta\neq 0$. Hence $\mult_{q^k}(C^k)=m$, and $C^k$ is transverse to both  $E_k$ and $L^k$ at  $q^k$.   
    \end{proof}

\begin{remark}
 If $k=1$  the first part of the statement is vacuous and the only interesting part is \eqref{lm-ka}, which states that
 $q^2\in C^2$ is an $(n-m,m)$-point with respect to the exceptional divisor $E_2$
\end{remark}

If $q\in C\subset \PP^2$  is a singular unibranch point, hence an $(m,n)$-point, one can  resolve the singularity by    a finite sequence of blow-ups of $\PP^2$, each centered at a point lying over $q$, such as the one described before Lemma~\ref{lm-k}. This process   
leads  to consider the euclidean algorithm for the pair $(m,n)$.
This connection is known; see \cite{Wall}, for example.
Our setting here is different, as  we need to extend our analysis to   smooth points of curves and to   tangent lines.

To a  pair of integers, $(m,n)$, with $n>m>0$, the euclidean algorithm associates a chain of integers determining their greatest common divisor, $c:=\gcd(m,n)$, as follows
\begin{equation}
 \label{eq-eu}
 l_0=n> l_1=m>l_2>\ldots >l_{r-1}>l_{r}= c >l_{r+1}=0
\end{equation}
where    
 $l_j$ is defined inductively for $ 2\leq j\leq r+1$  $$
  l_j:=l_{j-2}-l_{j-1}\lfloor l_{j-2}/l_{j-1}\rfloor.$$ 
We refer to \eqref{eq-eu} as the {\it euclidean sequence} of the pair $(m,n)$.
We set  
\[k_j=\lfloor l_{j-1}/l_j\rfloor, \quad    j=1,\dots, r.\]
As we shall see in Proposition~\ref{hypres}, in   the resolution sequence  resolving the $(m,n)$-point $q$, the first multiplicities
that appear over $q$  are  $l_1,l_2,\ldots, l_r$, in this order.
Moreover, each multiplicity $l_j$ appears exactly $k_j$ consecutive times.

We want to use this to give a  more convenient notation for the resolution sequence.
We define a total order on the following set  of pairs
  associated to the euclidean sequence  \eqref{eq-eu}
$$
\Pi_{(m,n)}:=\{(j, i),\quad \forall j=1,\ldots, r,\quad i=1,\ldots, k_j\}.
$$
We will use the index $i$ even though its range depends on $j$ to ease the notation.
We order $\Pi_{(m,n)}$ lexicographically as follows 
$$
(j',i')\geq (j,i)\  \text { if }\ j'>j \text { or } j'=j, \  i'\geq i.
$$
If $(j',i')$ is the minimum in $\Pi_{(m,n)}$ such that $(j',i')>(j,i)$, we say that $(j',i')$ is {\it next} to $(j,i)$, and write $(j',i')=\operatorname{next}(j,i)$. Explicitly:
$$
\operatorname{next}(j,i)=\left\{ \begin{array}{ll}
 (j,i+1)  &  \  \text {if  } i<k_j\\
 (j+1,1)  &  \  \text {if  } i= k_j. 
 \end{array} \right.
 $$
 We set
 $S^1_1:=\PP^2,\quad C^1_1:=C,\quad q^1_1:=q.$ 
 Now, $q^1_1$ is the center of the first blow-up
 $ 
 \sigma^i_j:S^i_j\to S^1_1
 $ 
 with $(j,i)=\operatorname{next}(1,1)$. We denote by $C^i_j\subset S^i_j$ the proper transform of $C$,
 and by $q^i_j\in C^i_j$ the point lying over $q$.
 We generalize to all $(j,i)>(1,1)$. 
We denote 
 $$q^i_j\in C^i_j \subset S^i_j\stackrel{\sigma^i_j}{\la} S^{i'}_{j'}\la \PP^2$$ 
so that $(j,i)=\operatorname{next}(j',i')$, the proper transform of $C$ is $C^i_j$, and $q^i_j$ is the only point of $C^i_j$ lying over $q$. The map $\sigma^i_j$ is the blow up at $q^{i'}_{j'}$, and we denote by $E^i_j \subset S^i_j$ its exceptional divisor.
For convenience, we denote $L=E^1_1\subset S^1_1=\PP^2$.
The divisors $E^1_j$ and their proper transform play a special role, hence 
for every $j=1,\ldots, r$ we denote by 
$$
 F^i_j\subset  S^i_j 
$$
 the proper transform of $E^1_j$, so that  $F^1_j=E^1_j$.
 We abuse notation and denote by  
$ 
\sigma^i_j:C^i_j \to  C^{i'}_{j'} 
$ 
the restriction of $\sigma^i_j$ to $C^i_j$.
 \begin{prop}
\label{hypres} 
Let    $q\in C\subset \PP^2$ be an    $(m,n)$-point,
  \eqref{eq-eu}   the euclidean sequence of $(m,n)$,  and   
 $\nu_{q}:C^{\nu}_q\to C$  the desingularization   of $C$ at $q$. 
We have a
  chain $\beta$ of   birational morphisms 
   \begin{align*}
  \beta:  C^{k_{r}}_{r}  \stackrel{\sigma^{k_r}_{{r}}}{\la}      \ldots  \to   C^1_{r}   \stackrel{\sigma^1_{{r}}}{\la}  C^{k_{r-1}}_{r-1}  \to&   \ldots    \quad  \quad  \quad\quad  \quad \quad  \quad  \quad\quad  \quad \quad  \quad \\
\ldots  \to  C^{k_j}_{j} \stackrel{\sigma^{k_j}_{j}}{\la}   \ldots \to C^{i}_{j} \stackrel{\sigma^{i}_{ j}}{\la}  C^{i-1}_{j}  \to 
\ldots \to  C^1_{j} \stackrel{\sigma^1_{ j}}{\la} C^{k_{j-1}}_{j-1}\to &  \ldots  \to   C^1_{1}=C 
  \end{align*}
  such that the following holds.
  \begin{enumerate}[(a)]
 \item
 \label{hypresa}
 For every $j=1,\ldots, r$ and  $i=1,\ldots k_j$ with $(j,i)\neq (r,k_{r})$
the point $q^i_j\in C^i_{j}$ is  an $(l_j,l_{j-1}-(i-1)l_j)$-point   w.r.t.   $F^i_j$.
\item
 \label{hypresb}
  $C^{k_r}_{r}$ has a $c$-fold point in $q^{k_r}_{r}$ where  it is transverse to   $E_{r}^{k_r}$ and to  $F_{r}^{k_r}$.
\item
 \label{hypresc}
If $c=1$ then $C^{\nu}_q=C^1_{r}$,  
and $q^1_{r}$ is a  $(1, l_{r-1})$-point w.r.t. $F^1_r$.
\item
 \label{hypresd}
 If $c>1$     then $\nu_{q}$ factors as follows  $\nu_{q}:C^{\nu}_q\la C^{k_r}_{r}\stackrel{\beta}{\la} C$.

\end{enumerate}
    \end{prop}
\begin{proof}
The chain was  defined before the statement. 
Of course, $q$ is an $(m,n)$-point w.r.t. $L=E^1_1$.
 If $n/m\in \N$ then $r=1$ and $c=m$, 
 and   the statement follows from  Lemma~\ref{lm-k}, which is indeed a special case of this Proposition.
 
Assume $n/m\not\in \N$, i.e. $k_1\neq n/m$.
 By Lemma~\ref{lm-k},  for every $i=1,\ldots, k_1$   the point  $q^i_1\in C_1^i$ is
 an $(m,n-(i-1)m)$-point   w.r.t. the proper transform,  
 $F_1^i$, of $E^1_1$; moreover,  $q^1_2\in C^1_2$ is an  $(l_2,l_1)$-point with respect to
  $E^1_2$. 
Now we continue inductively on $j$; assume   the statement for the level $j-1$ and  consider
  the $j$-th level of the chain, with $2\leq j\leq r-1$.
Then  $q^1_j\in C^1_{j}$  is  an 
 $(l_j,l_{j-1})$-point w.r.t.
  $E^1_j=F^1_j$, and $E^1_j$ is a line locally at $q^1_j$.  
By  the same lemma, for every $i=1,\ldots, k_j$   the point  $q^i_j\in C_j^i$ is
 an $(l_j,l_{j-1}-(i-1)l_j)$-point   with respect to the proper transform,  
 $F_j^i$, of $E^1_j$.  

Consider the last level,  $j=r$. Now $l_{r}=c$ and  $l_{r-1}/l_r=k_r$.
By Lemma~\ref{lm-k}, if  $i=1,\ldots ,k_r-1$ then 
$q^i_r\in C_r^i$ is
 an $(c,l_{r-1}-(i-1)c)$-point   w.r.t. the proper transform,  
 $F_r^i$, of $E^1_r$.   Moreover, by part \eqref{lm-kb} of the same lemma, 
$C_r^{k_r}$ has   a  $c$-fold point  in $q_r^{k_r}$,   and it is transverse to   $F_{r}^{k_r}$  and to the exceptional divisor.    The rest of the statement  is clear.
\end{proof}
  
We highlight  the following 
\begin{remark}
\label{rk-tan}
 For every $j\geq 1$  the curve $C^1_j$ is tangent to  $E^1_j$,  whereas  
$C^i_j$  meets  $E^i_j$  transversally for all $i\neq 1$.
 \end{remark}
   \subsection{Normal crossings resolutions}\label{subsec:ncr}
  Let $q\in C\subset \PP^2$ be an $(m,n)$-point and consider the sequence of maps in Proposition~\ref{hypres}. 
We    now add   the  blow-up  of   the last surface, $S_{r}^{k_r}$, at    
  $q_{r}^{k_r}$, written 
  $$\sigma_*:S_*\la S_{r}^{k_r}.$$    Its exceptional divisor will be written $E_*$; let  $C_*\subset S_*$ be the proper transform of $C$ and $q_*\in C_*$ the point lying over $q$.  We set $\beta_*= \beta\circ \sigma_*$ so that we have, abusing notation, 
$$
\beta_*:C_*\stackrel{\sigma_*}{\la}   C_{r} ^{k_r}\la  C_{r} ^{k_{r}-1}\la\ldots   \la  C^1_{r}\la \ldots   \la  C.
$$
 The subscript ``$_*$" indicates that, if  $c=\gcd(m,n)>1$, the type of singularity  of  $q_*$ is not known a priori  (we  know it is a $(c,h)$-point form some $h$)   in which case, by part \eqref{hypresd} of Proposition~\ref{hypres},  $\beta_*$ is only a partial resolution of $q$.
 \begin{example}
 \label{(3,5)}
Consider a curve $C\subset \PP^2$ with a $(3,5)$ point, $q$, with respect to the line $L$. The euclidean sequence is 
$5>3>2>1$ hence  $r=3$, and $k_1=k_2=1$, and $k_3=2$.
The picture below represents   the preimage of $C\cup L$ under four consecutive blow-ups. The first three blow-ups are those described in Proposition~\ref{hypres}, the last, introduced above,   is such that the resulting curve has normal crossings singularities. 
The first blow-up gives a $(2,3)$-point whose tangent   is the proper transform  of $L$.
The second blow-up gives a $(1,2)$-point  whose tangent  is the exceptional divisor. 

For better clarity,  the curve $C$ and its proper transforms are   curvilinear and red. The line $L$ and its proper transforms are horizontal  and  blue. 
The exceptional divisor of each blow-up is a vertical segment, and its proper transforms are segments (neither vertical nor horizontal).

 \begin{center}

 \begin{tikzpicture}[xscale=1.2,yscale=1.2]
  \node at (0,-1.5){$S_*$};
     \draw  [black] (0,1)--(0,-1)node at (0,1.2){$E_*$};
 \draw  [thick,blue] (-1.6,0.8)--(-0.3,0.8);
 \draw  [black] (-1,0.9)--(0.2 ,0.5);
  \draw[very thick, red] (-1,0.3)to  [out= 0, in= 120] (0.2,-0.3);
   \draw  [black] (-1,-1)--(0.2,-0.5);
  \draw  [black] (-1.5,-0.6)--(-0.2 ,-0.9);
  \node at (2,-1.5){ $S_3^2$};
 \draw  [thick,blue] (1,0.6)--(2.2,0.6);
  \draw[very thick, red] (1,0.3)to  [out= 0, in= 120] (2.2,-0.3);
   \draw  [black] (1,-0.8)--(2.2,0.1);
  \draw  [black] (1,-0.3)--(1.8 ,-0.8);
   \draw  [black] (2,1)--(2,-1)node at (2,1.2){$E_3^2$};
   \node at (4,-1.5){$S_3^1$};
   \draw  [thick,blue] (3,0)--(4.2,0);
  \draw[very thick, red] (4,0)to  [out=   100, in= 0] (3,0.7);
 \draw[very thick, red] (4,0)to  [out=   260, in= 0] (3,-0.7);
    \draw  [black] (3,1)--(4.2,0.5);
   \draw  [black] (4,1)--(4,-1)node at (4,1.2){$E_3^1$};
      \node at (6,-1.5){$S_2^1$};
  \draw  [thick,blue] (5,0)--(6.2,0);
  \draw[very thick, red] (6,0)to  [out=   170, in= 290] (5,0.7);
 \draw[very thick, red] (6,0)to  [out=   190, in= 70] (5,-0.7);
  \draw  [black] (6,1)--(6,-1)node at (6,1.2){$E_2^1$};
  \node at (8,-1.5){$ \PP^2$};
 \draw  [thick, blue] (7.2,0)--(8.5,0)node at (7.2,0.2){$ L$};
  \draw[very thick, red] (8.5,0)to  [out=   180, in= 270] (7.2,0.7)  node at (8.3,0.3){$ C$}; 
  \node at (8.4,-0.2){$q$};
 \draw[very thick, red] (8.5,0)to  [out=   180, in= 90] (7.2,-0.7);
     \end{tikzpicture}
 \end{center}
 This example is a   Fibonacci point, described in general  in example~\ref{Fibonacci}.
  \end{example}

We write  
 \begin{equation}
\label{eq-D*}
 D_*=D_*(q):=(\beta_*^{-1}(C\cup L))_{\operatorname{red}}\subset S_*
 \end{equation}  for  
the reduced subscheme/divisor supporting $\beta_*^{-1}(C\cup L)$. In the above picture, $D_*$ is the first   from the left.

 We write
 $A^i_{j}  \subset S_*$ for the proper transform of $E^i_{j}$ in $S_*$, so that  $A^1_1$ is the proper transform of $L=E^1_1$. Then  
\begin{equation}
 \label{def-D*}
 D_*= C_*+E_* +\sum_{\stackrel{j=1,\ldots, r}{i=1,\ldots, k_j}}  A^i_{j}.
\end{equation}
By construction,  two components of $D_*$ intersect   in at most   one point, and no three of them intersect in the same point.
 
\begin{lemma}
 \label{lm-q*} Notation as above. Then

\begin{enumerate}[(a)]
 \item \label{lm-q*a}
 $C_*\cap A^i_j=\emptyset$ for  every $(j,i)$.
 \item\label{lm-q*b}
 $E_*\cap A^i_j\neq \emptyset$ if and only if  $(j,i)\in\{ (r,k_r),(r,1)\}.$

\item\label{lm-q*c}
For every $(j,i)\not\in\{ (r,k_r),(r,1)\}$ there is a unique $(j',i')>(j,i)$ such that $A^i_j\cap A^{i'}_{j'}\neq \emptyset$, moreover
 \begin{equation}
 \label{eq-ij}
(j',i')=\left\{ \begin{array}{ll}
\operatorname{next}(j,i)  &  \  \text {if  } i\neq 1 \\
\operatorname{next}(j+1,1)  &  \  \text {if  } i= 1. 
 \end{array} \right.
\end{equation}
 \item
 \label{lm-nc}
 The divisor $D_*-C_*$ is normal crossings.
If $c=1$  and $C$ is normal crossings away from $q$, then   the divisor 
  $D_*$ is   normal crossings.   
\end{enumerate}
 \end{lemma}
 \begin{proof} 
 By Proposition~\ref{hypres},
the  curve $C_r^{k_r}$ has a unibranch $c$-fold point in $q_{r}^{k_r}$, where it meets transversally both the proper transform of $E^1_r$ and the exceptional divisor, $E_r^{k_r}$. Hence the exceptional divisor, $E_*$, of  the  blow-up in $q_{r}^{k_r}$, intersects $C_*$ and     the proper transforms of $E^1_r$ and of $E_r^{k_r}$ (i.e. $A_{r}^{1}$  and  $A_r^{k_r}$) in three different points.
Hence the first two statements follow,
 
To prove \eqref{lm-q*c}, let  $(j',i')>(j,i)$; notice that $A^i_j\cap A^{i'}_{j'}\neq \emptyset$ if and only if for  some  $S=S^{i''}_{j''}$
with $(j'',i'')\geq (j',i')$    the proper transforms of $E^i_j$ and $E^{i'}_{j'}$ in $S$ intersect away from the proper transform of $C$.
Suppose $i\neq 1$;  Proposition~\ref{hypres}  implies that  $E^i_j$ is transverse to $C^i_j$, hence in the next blow-up, the proper transform of $E^i_j$ intersects the exceptional divisor away from the proper transform of $C$. 
Suppose $i=1$; then $E^1_j$ is tangent to $C^1_j$,  and    its proper transform is tangent to $C^i_j$ for every $i\leq k_j$.
In $S^1_{j+1}$ the curve $C^1_{j+1}$ is tangent to $E^1_{j+1}$, hence the proper transform of $E^1_j$ intersects
$E^1_{j+1}$ transversally in $q^1_{j+1}$. Therefore in the next blow-up, the proper transform of $E^1_j$ intersects the exceptional divisor away from the proper transform of $C$. \eqref{lm-q*c} is proved.

All components of $D_*-C_*$ are smooth,   intersect transversally, and only pairwise.
Hence $D_*-C_*$ is normal crossings.
If $c=1$, then    $q_{r}^{k_r}$ is a smooth point of $C_r^{k_r}$, hence $D_* $ is normal crossings if so is $C$ away from $q$. This proves (d) and concludes the proof.
 \end{proof}

 \section{Tropical curves of unibranch points}
 In this section we introduce the ``contact'' tropical curve associated to a unibranch point of a curve.
 \label{sec:tropical}
 \subsection{Dual graphs of partial resolutions}
  \label{sec-trop}
 Let $C\subset \PP^2$  have a unique (for simplicity)  singular point in $q$.  It is well known that there exists a unique surface $S_\rho$ with a birational morphism $\rho:S_\rho\to \PP^2$ such that the strict transform of $C$ in $S_\rho$ is smooth,  the reduced scheme underlying $\rho^{-1}(C)$ is a divisor with normal crossings,  written $D_\rho$, and $S_\rho$ is minimal with respect to these properties. The dual graph, $G_\rho$, of $D_\rho$
 (whose vertices   are the irreducible components of $D_\rho$, with 
  an edge between two vertices if the two components intersect)
 is often called the resolution graph of $q$;  \cite{Wall}.
The graph  $G_\rho$ is not   suitable here, and we will   define a variant of it which will lead, in turn,  to   a tropical curve.
 
 We    consider the dual graph, $$G_*=G(D_*(q))$$  of the divisor   $D_*$, defined in Section~\ref{subsec:ncr}.
  The vertices of $G_*$  are the irreducible components of $D_*$, hence $G_*$ has $2+ \sum_{i=1}^rk_i$ vertices,
 and there is an edge between two vertices if the two components intersect. We shall abuse notation and denote the vertices and the components they represent by the same symbol.
 For the next lemma,  recall that a {\it leaf} of a graph is a vertex of degree   $1$ (and the degree, or valency,  of a vertex is the number of half-edges adjacent to it).

\begin{lemma}
 \label{lm-G}
Let $q\in C$ be an $(m,n)$-point.
\begin{enumerate}[(a)]
\item
 \label{lm-G2}
The graph   $G(D_*(q))$ is a tree with    a unique vertex, $E_*$,  of degree $3$
and 
 three leaves among which $C_*$ and $A^1_1$, and has the following form
 $$
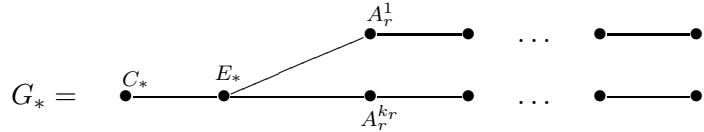
\xymatrix@=1.pc{
&&&&&&*{\bullet} \ar@{-}[rr]^(.1){A^1_r}  &&*{\bullet}& \ldots  &    *{\bullet}\ar@{-}[rr] &&*{\bullet}&&\\
G_*=&*{\bullet} \ar@{-}[rr]^(.1){C_*} &&*{\bullet} \ar@{-}[rrru]^(.1){E_*}\ar@{-}[rrr] &&&
*{\bullet} \ar@{-}[rr]_(.1){A^{k_r}_r} &&*{\bullet}  &\ldots    & *{\bullet}\ar@{-}[rr] &&*{\bullet}\\
&&&&&&&&&&&
}
$$
\item
 \label{lm-GG}
 $G(D_*(q))$ only depends on the pair $(m,n)$, i.e.
if $q'\in C'\subset S'$ is an $(m,n)$-point, then 
 $G(D_*(q))$  is isomorphic to $G(D_*(q'))$ via an isomorphism that takes the vertex $C_*$ to the vertex $C'_*$.
\end{enumerate} \end{lemma}

\begin{proof}
Write $G_*=G(D_*(q))$.
 (The fact that $G_*$ is a tree should to   be clear from its definition, but we give it a proof   below.) 
 
Lemma~\ref{lm-q*} implies that $E_*$ is the only vertex of degree $3$, and it is joined by an edge to $C_*$,  $A_r^{k_r}$ and $A^1_r$, as shown in the picture. The same lemma implies that $C_*$ and  $A^1_1$ are leaves of $G_*$ (and that  the third leaf is $A^1_2$ if $k_1=1$, and $A^2_1$ if $k_1>1$).

Let $\nu_i$ be the number of vertices of degree $i$ in $G_*$. Then $\nu_3=1$ and $\nu_1\geq 2$. The number of edges of $G_*$ is thus
 $$
 |E(G_*)|=(\nu_1+2\nu_2+3\nu_3)/2=(\nu_1+2\nu_2+3)/2
 $$
  hence $\nu_1$ must be odd, hence $\nu_1\geq 3$.
Since $G_*$ is connected, its genus $g(G_*) = b_1(G_*)$ satisfies $ g(G_*)= |E(G_*)|-  |V(G_*)|+1$. Moreover, since $|V(G_*)|=\nu_1+\nu_2+1$, we have 
  $$
  g(G_*) =(\nu_1+2\nu_2+3)/2-(\nu_1+ \nu_2+1)+1=(3-\nu_1)/2\leq 0
  $$
  with equality only if $\nu_1=3$.
Since $g(G_*)\geq 0$ we get $\nu_1=3$ and $g(G_*)=0$.
 \eqref{lm-G2} is proved.
Part \eqref{lm-GG} is clear.
 \end{proof}

As  $G(D_*(q))$ depends only on the pair $(m,n)$,   we set
$
 G^{(m,n)}_*:=G(D_*(q)) 
$
 for any $(m,n)$-point $q$ on any curve $C$. We will show later (see Lemma \ref{lm-mnequiv}) that it depends only on the euclidean sequence of $(m,n)$.
 
Different pairs  can give the same graph, as in the next example.
 
\begin{example}
\label{ex-path}
Let  $q$ be a smooth point of type $(1,4)$, hence  
 $r=1$ and  $k_1=4$.
The graph
 $G(D_*(q))$ is as follows
 $$
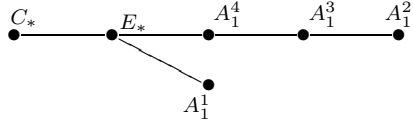
\xymatrix@=.8pc{
&*{\bullet} \ar@{-}[rr]^(.1){C_*} &&*{\bullet} \ar@{-}[rr]  \ar@{-}[rrd]^(.1){E_*}_(1){A^1_1}  &&*{\bullet} \ar@{-}[rr]^(.2){A^{4}_1}  &&*{\bullet}\ar@{-}[rr]^(.2){A^{3}_1} ^(1){A^{2}_1} && *{\bullet}   \\
&&&&&*{\bullet}&&&&&&&&&&&&&&&&&&&&&&&&&&&&
}
$$
Let now   $q'$ be  a $(2,8)$-point,  then $q'$ is   singular, $r=1$ and $k_1=4$,
 and the  graph $G(D_*(q'))$  is the same as above.   \end{example}

\begin{example}{\it Fibonacci singularities.}
\label{Fibonacci}
A {\it Fibonacci point} is an $(m,n)$-point   such that   $(m,n)$ is a pair of consecutive integers in the Fibonacci sequence $\{f_n\}$ defined by $f_1=1$, $f_2=2$, and $f_n=f_{n-1}+f_{n-2}$:
$$
 1,\ 2,\ 3,\  5,\ 8,\  13\ , 21, \ \ldots 
$$
The ``blow-up''  of a Fibonacci point  of type $(f_n, f_{n-1})$ is   a Fibonacci point of type $(f_{n-1}, f_{n-2})$.
The euclidean sequence for a Fibonacci $(m,n)$-point is equal to the reversed Fibonacci sequence truncated at $n$, for example for
  $(m,n)=(13,21)$ we get
$$
l_0=21>l_1=13>l_2=8>l_3=5>l_4=3>l_5=2>l_6=1.
$$
An $(m,n)$-point is a Fibonacci point if and only if $k_r=2$ and $k_j=1$ for all $j=1,\ldots, r-1$.
For such a  point  the chain of blow-ups   in \ref{hypres} is
$$
C_*\la C^2_r\la  C^1_r\la C^1_{r-1}\la\ldots \la C^1_2\la C.
$$
The   graph   $G_*$ for Example~\ref{(3,5)} is the following
$$\xymatrix@=.8pc{
&&&&&*{\bullet} \ar@{-}[rr]^(.1){A^1_3} ^(1){A^1_1} &&*{\bullet} \\
&*{\bullet} \ar@{-}[rr]^(.1){C_*} &&*{\bullet} \ar@{-}[rru]^(.1){E_*}\ar@{-}[rrd] &&&&   \\
&&&&&*{\bullet} \ar@{-}[rr]_(.1){A^2_3} _(1) {A^1_2}&&*{\bullet} &&\\
&&&&&&&&&&&\\
}
$$
If $r$ is even, $G_*$ is as follows
 $$\xymatrix@=.8pc{
&&&&&*{\bullet} \ar@{-}[rr]^(.1){A^1_r}  &&*{\bullet} \ar@{-}[rr]^(.2){A^1_{r-2}}  &&*{\bullet}&\ldots   & *{\bullet}\ar@{-}[rr]^(1){A^1_2}&&*{\bullet}&&\\
&*{\bullet} \ar@{-}[rr]^(.1){C_*} &&*{\bullet} \ar@{-}[rru]^(.1){E_*}\ar@{-}[rrd] &&&&   \\
&&&&&*{\bullet} \ar@{-}[rr]_(.1){A^2_r} &&*{\bullet}  \ar@{-}[rr]_(.2){A^1_{r-1}}&&*{\bullet}  \ar@{-}[rr]_(.2){A^1_{r+1}}&&*{\bullet} &\ldots & *{\bullet}\ar@{-}[rr]_(1){A^1_1}_(.1){A^1_3}&&*{\bullet}\\
&&&&&&&&&&&\\
}
$$
  \end{example} 

We now study under what conditions different pairs give the same graph.
\begin{defi}
\label{def-eq} 
Let $(m,n)$ and $(m',n')$ be two pairs of integers such that $n>m>0$ and $n'>m'>0$.
The euclidean sequence for $(m,n)$ is given in \eqref{eq-eu};    let the euclidean sequence of $(m',n')$ be
\begin{equation}
 \label{eq-eu'} n'= l'_0> m'=l'_1>\ldots l'_{r}=c'>l'_{r'+1}=0\end{equation}
    where $c'=\gcd(m',n')$. Set  $k'_j =\lfloor l_{j-1}/l_j\rfloor$ for  $j=1,\dots, r'.$
   We  say that the two pairs are {\it equivalent}, and write
    $$(m,n)\sim (m',n'),$$
    if $r=r'$ and $k_j=k'_j$ for every $j=1,\dots, r.$ 
    \end{defi}
   If $q\in C\subset \PP^2$ is an $(m,n)$-point and  $q'\in C'\subset \PP^2$ is an $(m',n')$-point we say that $q$ and $q'$ are equivalent if $(m,n)\sim (m',n')$.
 
\begin{lemma}
\label{lm-mnequiv}
 If $(m,n)\sim (m',n')$ then $G^{(m,n)}_*\cong G^{(m',n')}_*$.
\end{lemma}
\begin{proof}
 By Proposition~\ref{hypres}, for any $(m,n)$-point $q$,  the divisor $D_*(q)$ is completely determined by the integers $r$ and $k_1,\ldots, k_r$ of the euclidean sequence of $(m,n)$, hence $G^{(m,n)}_*$ depends only on the equivalence class of $(m,n)$.
\end{proof}
  The converse fails, as the following example shows.
\begin{example}
\label{ex-2g}
If $(m,n)=(5,8)$, then $q$ is a Fibonacci point with $r=4$ and the graph   is as follows
$$
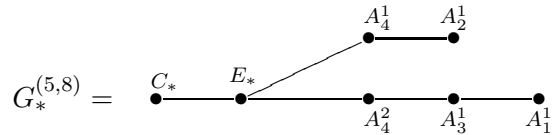
\xymatrix@=.8pc{
&&&&&&*{\bullet} \ar@{-}[rr]^(.1){A^1_4}^(1){A^1_2}  &&*{\bullet}  && &&&&\\
G^{(5,8)} _*=&*{\bullet} \ar@{-}[rr]^(.1){C_*} &&*{\bullet} \ar@{-}[rrru]^(.1){E_*}\ar@{-}[rrr] &&&*{\bullet} \ar@{-}[rr]_(.1){A^{2}_4}_(1){A^1_3} &&*{\bullet}\ar@{-}[rr]_(1){A^1_1}    &&*{\bullet}&&\\
&&&&&&&&&&&\\
}
$$
Now let $(m',n')=(3,7)$, then $r'=2$ hence $(m,n)\not\sim (m',n')$. The   graph   below is clearly isomorphic to $G^{(5,8)}_*$
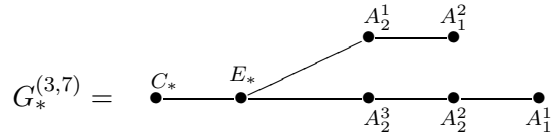
$$\xymatrix@=.8pc{
&&&&&&*{\bullet} \ar@{-}[rr]^(.1){A^1_2}^(1){A^2_1}  &&*{\bullet}  && &&&&\\
G^{(3,7)}_* =&*{\bullet} \ar@{-}[rr]^(.1){C_*} &&*{\bullet} \ar@{-}[rrru]^(.1){E_*}\ar@{-}[rrr] &&& *{\bullet} \ar@{-}[rr]_(.1){A^{3}_2}_(1){A^2_2} &&*{\bullet} \ar@{-}[rr]_(1){A^1_1}    &&*{\bullet}   &&&&\\
&&&&&&&&&&&\\
}
$$
\end{example}
It is well known that the condition  $(m,n) \sim (m',n')$   can be expressed more directly as follows.
    \begin{remark}
\label{rm-mnfrac}
  $(m,n)\sim (m',n')$ if and only if  $m/n=m'/n'$.
\end{remark}
See \cite[Subsection 10.6 and Theorem 160]{HW}, for example.
      \subsection{The contact tropical curve of an $(m,n)$-point}
 To improve upon Lemma~\ref{lm-mnequiv} and give a geometric characterization for equivalence of pairs, we now introduce a tropical curve associated to an $(m,n)$-point.
 We first put 
on the graph $G_*=G(D_*(q))$   a length   on its edges, 
 $\ell_*:E(G_*)\to \R_+$,
such that $\ell_*(e)=1$ for every edge $e$ of  $G_*$.
This defines  a tropical curve $\Gamma_*:=(G_*,\ell_*)$.
We   shall now define the {\it contact tropical curve} of the $(m,n)$-point $q$, written  
 $\Gamma_q$. The word ``contact" is motivated by Theorem~\ref{mirror}. 
  As for notation, set 
$$\Gamma_q=(G_q,\ell_q ),\quad \quad G_q=(V(G_q),E(G_q), L(G_q))$$ 
with $G_q$ a graph whose vertex set, edge set, and leg set are  written, respectively, $V(G_q)$, $E(G_q)$ and  $L(G_q)$
(a leg   is a half-edge attached to only one vertex),
and $\ell_q$ is a   function
$$
\ell_q: E(G_q) \cup  L(G_q)\la \R_+.
$$
We define  $\Gamma_q$  as the tropical curve    
 obtained from $\Gamma_*$ by removing every vertex  $A^i_j$ with $i>1$
 without disconnecting the graph and adding the corresponding lengths. 
 Since the edges   $A^i_j$ have degree at most $2$, this operation is well defined. Indeed, if 
$A^i_j$  has degree $1$, then the
  edge adjacent to it  becomes a leg, or part of a leg, of $\Gamma_q$;  if $A^i_j$ has degree $2$  the
  two    edges  adjacent to $A^i_j$ are merged in a new edge or   in a leg.

\begin{remark}
 By Remark~\ref{rk-tan}, we  removed  from $G_*$ those vertices corresponding to exceptional divisors  $E^i_j$ that are not tangent to the proper transform $C^i_j$ of $C$.
\end{remark}
  
We have  an obvious surjection  from the edges of $G_*$ to the  edges and legs  of $G_q$, written
  $$
\mu:E(G_*)\la E(G_q)\cup L(G_q).
  $$
The length function $\ell_q$ is defined as follows  $$
\ell_q(x):=\sum _{e\in \mu^{-1}(x)}\ell_*(e)
$$
  for any $x\in E(G_q)\cup L(G_q)$.
  The  definition of the tropical curve $\Gamma_q$ is now complete. Let us describe it more closely.
Its
  vertex set  is  
$$ V(G_q)=\{C_*,\  E_*,\  A^1_j,\quad   j=1,\ldots,r\}.$$
An edge of $G_q$ can be written   unambiguously  as
  $e=vw$ with $v,w\in V(G_q)$, hence, by Lemma~\ref{lm-q*}, we have
$$ E(G_q)= \left\{ \begin{array}{ll}
 \{C_* E_*,\  E_*A^1_r,\ E_*A^1_{r-1},\  A^1_{j+1} A^1_{j-1},\  \forall j=2,\ldots, r-1\}   &  \   \text {if  } r>1\\
 \{C_* E_*,\  E_*A^1_1\}  & \   \text {if  } r=1\\
 \end{array} \right.$$
 and one easily checks the following
  \begin{equation}
 \label{eq-lGq}
\ell_q(vw):=\left\{ \begin{array}{ll}
 1  &  \   \text {if  } vw= E_*C_*  \text{ or } vw= E_* A^1_r\\
 k_r & \   \text {if  } vw=  E_*A^1_{r-1}\\
k_j  &  \  \text {if  } vw= A^1_{j+1} A^1_{j-1},\quad \forall j=2,\ldots, r-2.\\
 \end{array} \right.
\end{equation}
Finally,  if $r>1$ and $k_1>1$ then $G_q$ has a leg attached to $A^1_2$ of length $k_1-1$. If $r=1$ then $G_q$ has a leg attached to $E_*$ of length $k_r-1$. 
In particular,  $\Gamma_q$ has at   most one leg.
Summarizing and applying  Lemma~\ref{lm-G}:
 \begin{lemma}
 \label{lm-Gamma}
Let $q\in C\subset S$ be an $(m,n)$-point and $\Gamma_q=(G_q,\ell_q)$ its contact tropical curve.

\begin{enumerate}[(a)]
\item
 $G_q$ is a tree with  $r+2$ vertices, and  a unique vertex, $E_*$,  of degree $3$. 
 
 The vertices $C_*$ and $A^1_1$ are   leaves of $G_q$ (with no leg attached).

 $G_q$ has no legs if $k_1=1$, and a unique leg of length $k_1-1$ if $k_1>1$.

 
\item If $q'\in C'\subset S'$  is an $(m,n)$-point, then  the contact tropical curves $\Gamma_q$  and  $\Gamma_{q'}$ are isomorphic.
 \end{enumerate} \end{lemma}
 As    $\Gamma_q$ only depends on  $(m,n)$  we set, for any $(m,n)$-point $q\in C$,
$$
 \Gamma^{(m,n)}:=\Gamma_q.
$$

 \begin{example}
\label{ex-2t}
Here are the tropical curves corresponding to   Example~\ref{ex-2g}.
$$\xymatrix@=.8pc{
&&&&&*{\bullet} \ar@{-}[rr]^(.1){A^1_4}^(1){A^1_2}^{\color{blue}1}  &&*{\bullet}  && &&&&\\
\Gamma^{(5,8)}   =&*{\bullet} \ar@{-}[rr]^(.1){C_*}_{\color{blue}1} &&*{\bullet} \ar@{-}[rru]^(.1){E_*}^{\color{blue}1}\ar@{-}[rrrr]_{\color{blue}2}_(1){A^1_3}  
 && &&*{\bullet}  \ar@{-}[rr]_(1){A^1_1}_{\color{blue}1}&& *{\bullet} &&\\
}$$
and
 $$\xymatrix@=.8pc{
&&&&&*{\bullet} \ar@{-}[rr]^(.1){A^1_2}^{\color{blue}1}   &&  && &&&&\\
\Gamma^{(3,7)}  =&*{\bullet} \ar@{-}[rr]^(.1){C_*}_{\color{blue}1} &&*{\bullet} \ar@{-}[rru]^(.1){E_*}^{\color{blue}1}\ar@{-}[rrrrr]_{\color{blue}3}_(1){A^1_1} 
 &&&&& *{\bullet}   &&&&\\
}
$$
These tropical curves are not isomorphic, while   $G^{(5,8)}_*\cong G^{(3,7)}_*$.
\end{example}
 
 We can now geometrically characterize equivalent pairs.
\begin{thm}
\label{thm-mn}
  $(m,n)\sim (m',n')$ if and only if  $\Gamma^{(m,n)}\cong \Gamma^{(m',n')}$ (i.e. $\Gamma^{(m,n)}$ and $\Gamma^{(m',n')}$ are isomorphic as tropical curves).
\end{thm}
\begin{proof}
 If $(m,n)\sim (m',n')$ then Lemma~\ref{lm-mnequiv} yields $G^{(m,n)}_*\cong G^{(m',n')}_*$, hence, by what we saw above,        $\Gamma^{(m,n)}\cong \Gamma^{(m',n')}$.
 
To prove the converse, write $\Gamma=\Gamma^{(m,n)}$ and $\Gamma'=\Gamma^{(m',n')}$ and  assume  that there is an isomorphism of tropical curves $\phi: \Gamma \to  \Gamma'$. Hence  $\phi$ is a length-preserving isomorphism between the underlying graphs.

 Therefore  $\Gamma$ and $\Gamma'$ have the same number of vertices,  equal to $r+2$, hence $r=r'$. Moreover,   $\Gamma$ and $\Gamma'$ have the same number of legs, equal to $\min\{1, k_1-1\} = \min\{1, k'_1 -1\}$;
 if this 
  number is $0$ then  $k_1=k'_1=1$. If  $\min\{1, k_1-1\}=1$   then $\Gamma$ and $\Gamma'$ have each  one leg of respective lengths $k_1-1$
  and  $k'_1-1$. Since the length of the leg is  preserved by the isomorphism $\phi$, we get $k_1=k'_1$.
 We  write
 $$
V(\Gamma')= \{C'_*,\  E'_*,\  (A^1_j)',\quad \forall j=1,\ldots,r\}.
 $$
  Recall that  $\Gamma$ and $\Gamma '$ have each a unique vertex, respectively   $E_*$ and $E'_*$, of degree $3$. Therefore the isomorphism $\phi$ has to map  the vertex $E_*$ to  the vertex $E'_*$.
 If $r=1$  then $k_1\geq 2$, hence $\Gamma$ has a unique leg and looks as follows
 $$
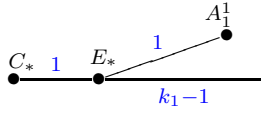
\xymatrix@=.8pc{
& &&&&&*{\bullet}    &&& \\
 &*{\bullet} \ar@{-}[rr]^(.1){C_*}^(.5){\color{blue}1} &&*{\bullet} \ar@{-}[rrru]^(1){ A^1_1} ^(.1){E_*}^(.5){\color{blue}1}\ar@{-}[rrrr]_(.5){\color{blue} k_1-1}&&&&&&&   \\
}
$$
Since   we already proved that   $k_1=k'_1$,     we are done.

If $r\geq  2$,  Lemma~\ref{lm-Gamma}  implies that    $\Gamma$ has the following form
  $$
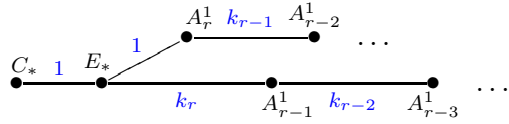
\xymatrix@=.8pc{
 &&&&*{\bullet} \ar@{-}[rrr]^(.1){ A^1_r}^(1){ A^1_{r-2}}^(.5){\color{blue} k_{r-1}}  &&&*{\bullet}&\ldots& \\
 *{\bullet} \ar@{-}[rr]^(.1){C_*}^(.5){\color{blue}1} &&*{\bullet} \ar@{-}[rru]^(.1){E_*}^(.5){\color{blue}1}\ar@{-}[rrrr]_(.5){\color{blue} k_r} 
&&&& *{\bullet}\ar@{-}[rrr]_(.1){A^1_{r-1}}_(1){A^1_{r-3}}_(.5){\color{blue} k_{r-2}}  &&&*{\bullet}&\ldots && \\
&&&&&&&&&&&\\
}
$$
and the same holds for $\Gamma'$, mutatis mutandis.
 Since $k_r\geq 2$,  the edges $E_*A^1_r$ and $E_*A^1_{r-1}$ have different lengths, equal to $1$ and $k_r$. The same holds on $\Gamma'$: as   $k'_r\geq 2$  the edges $E'_*(A^1_r)'$ and $E'_*(A^1_{r-1})'$ have   lengths  equal to $1$ and $k'_r$. Hence
 an isomorphism from  $\Gamma$ to $\Gamma'$  
 maps $A^1_r$ to $(A^1_r)'$ and $A^1_{r-1}$ to $(A^1_{r-1})'$, hence    $k_r=k'_r$. 
 Hence the edge $A^1_rA^1_{r-2}$ is mapped to the edge $(A^1_r)'(A^1_{r-2})'$, hence $k_{r-1}=k'_{r-1}$.
 Iterating,   the edge $A^1_jA^1_{j-2}$ is mapped to the edge $(A^1_j)'(A^1_{j-2})'$, hence $k_{j-1}=k'_{j-1}$
 for all $j=3,\ldots r$.  The proof is complete.
 \end{proof}

   \section{Hypertangency}
   \label{sec-mirror}
 \subsection{Hypertangency at unibranch point}
 Let  $B,C\subset S$ be  two  integral  curves and  $q\in B\cap C$;   let $\sigma:S'\to S$ be the blow up at $q$. If $B',C'\subset S'$  are   
  the proper transforms of $B$ and $C$,
    the following is well know   and will be used often  
 \begin{equation}
 \label{eq-int}
(B'\cdot C') =(B\cdot C) -\mult_q(B)\mult _q(C). 
\end{equation}
 
 For  our next result, we shall assume that $q$ is an $(m,n)$-point of $C$ and an $(m',n')$-point of $B$.
Recall, from Definition~\ref{def-eq}, that  two pairs $(m,n)$ and $(m',n')$ are   equivalent, written
 $(m,n)\sim (m',n'),$ 
    if $r=r'$ and $k_j=k'_j$ for every $j=1,\dots, r$, where the integers $r,r', k_j, k'_j$ are determined by the euclidean sequences of the   pairs. 
    
The proof of the next statement combines our earlier  results   with results  in \cite{GBP}.
  \begin{thm}
\label{mirror}
 Let $B, C\subset \PP^2$ be   integral curves 
 of  degree at least $2$. 
Let $q$ be an   $(m',n')$-point for $B$ and an  $(m,n)$-point for $C$  such that
 $B\cap C = \{ q\}.$ 
  Then the following equivalent conditions hold
  
\begin{enumerate}[(a)]
 \item
 $ \Gamma^{(m,n)}\cong \Gamma^{(m',n')}$.
 \item
$ m/n=m'/n'$.
 \item
 $(m',n')\sim (m,n)$.
  \end{enumerate}
   \end{thm}
 
 \begin{proof}
 The equivalence of the three conditions follows from   Remark~\ref{rm-mnfrac} and Theorem~\ref{thm-mn}.

 Set $b=\deg B$ and $d=\deg C$.
  By hypothesis  $(B\cdot C)_q=(B\cdot C)=bd. $ 
The assumption implies that $B$ and $C$ have the same tangent line, written $L$, at $q$, so that 
$(C\cdot L)_q=n$ and $(B\cdot L)_q=n'$.

To prove that $m/n=m'/n'$ we apply the  results in subsection 2.3 of  \cite{GBP},   in particular the ``Strong Triangle Inequality" in  Theorem 2.8 and Corollary 2.9. They imply that the smallest among the following three numbers  are equal
\begin{equation}
 \label{GBP}
\frac{(B\cdot C)_q}{\mult_q(B)\mult_q(C)},\quad\  \frac{(B\cdot L)_q}{\mult_q(B)\mult_q(L)}, \quad\  \frac{(C\cdot L)_q}{\mult_q(C)\mult_q(L)}.
\end{equation}
We have, 
$$
\frac{(B\cdot C)_q}{\mult_q(B)\mult_q(C)}=\frac{ bd}{m'm}\geq \frac{ bn}{m'm}> \frac{ m'n}{m'm}=\frac{n}{m}=\frac{(C\cdot L)_q}{\mult_q(C)\mult_q(L)}
$$
(using $d\geq n$ and $b>m'$  for the two inequalities).
Similarly
$$
\frac{(B\cdot C)_q}{\mult_q(B)\mult_q(C)} >  \frac{n'}{m'}=\frac{(L\cdot B)_q}{\mult_q(L)\mult_q(B)}.$$
Hence the second and third number in \eqref{GBP} are the smallest, hence  they are equal, hence $n/m=n'/m'$ and we are done.

\end{proof}

We now give a  different, self-contained proof which applies  the material of the previous sections.

\begin{proof}[Alternative proof]
We shall   prove   $(m',n')\sim (m,n)$.  We can assume $r\leq r'$.
Consider the chain of blow-ups described before  Proposition~\ref{hypres} for the curve $C$.
For every  $j=1,\ldots, r$ and $i=1,\ldots, k_j$  we have a sequence of  blow-ups, $S^i_j$, of $\PP^2$,   containing the proper transform, $C^i_j$, of $C$, and a point
  $q^i_j\in C^i_j$, which is the center of the blow up of $S^i_j$; we have    $q=q^1_1$,\  $C=C^1_ 1$ and $\PP^2=S^1_1$.  We    consider also the curve $C^{}_*$ in the blow-up  $S^{}_*$ of $S^{k_r}_r$ at $q^{k_r}_r$, defined before Lemma~\ref{lm-q*}. 
  
We denote by 
  $B^{(i)}_{(j)}\subset S^i_j$ and $B_{(*)}\subset S^{}_*$ the proper transforms of $B$.
  Since $C$ and $B$ meet only in $q$, if  the curves  $C_j^i$  and  $B^{(i)}_{(j)}$ intersect, they   do so in the only point of $C_j^i$ mapping to $q$, i.e. in $q_j^i$. 
Notice that if   $(m,n)\sim (m',n')$  
  the chain of blow-ups for $B$ described  in Proposition~\ref{hypres} corresponds exactly to the one  for  $C$, hence  $B^{(i)}_{(j)}=B^i_j$.

  We   prove    $(m,n)\sim (m',n')$  by contradiction.  
Then, setting 
 $$
h=\min\{j: k_j\neq k'_j,\  j=1,\ldots,r\},
$$
we have   $h\leq r$.
We  can  assume $k_h<k'_h$, the case $k_h>k'_h$ is similar. By the definition of $h$, we have  $B^{k_h}_h= B^{(k_h)}_{(h)}$ and
$C^{k_h}_h\cap B^{k_h}_h=\{q^{k_h}_h\}$.

\noindent
{\bf Claim.} {\it If  $h<r$  then $C^1_{h+1}$ is not tangent to 
$B^{(1)}_{(h+1)}$; if $h=r$ then 
$C^{k_r}_r$  is not tangent to $B^{(k_r)}_{(r)}$.}

Recall that, by   Proposition~\ref{hypres},  for every $(j,i)$ with $i>1$, except the pair $(r,k_r)$, the curve $C^i_j$ is tangent to $F^i_j$ (the proper transform in $S^i_j$
of the exceptional divisor of $S^1_j$), and   $C^1_j$ is tangent to the exceptional divisor $E^1_j\subset S^1_j$. 
If  $h< r$  both  $C^{k_h}_h$ and $B^{(k_h)}_{(h)}=B^{k_h}_h$    are tangent to $F^{k_h}_h$, hence their proper transforms in the 
  blow-up   at  $q^{k_h}_h$ intersect in $q^1_{h+1}$. Now, the proper transform of $C$ is $C^1_{h+1}$ and the proper transform of 
$B$, in the notation of Proposition~\ref{hypres}, is $B^{k_h+1}_h$ (as $k'_h>k_h$),   and we have  $B^{k_h+1}_h=B^{(1)}_{(h+1)}$.
In $q^1_{h+1}$, the curve $C^1_{h+1}$ is tangent to $E^1_{h+1}$, whereas $B^{k_h+1}_h$  is tangent to the proper transform of 
$F^{k_h}_h$, hence  $C^1_{h+1}$ is not tangent to 
$B^{(1)}_{(h+1)}$, as claimed.

If $h=r$,  
according to  Proposition~\ref{hypres}, the curve $C^{k_r}_r$ is not tangent to $F^{k_r}_r$ in  $q^{k_r}_r$, while $B^{k_r}_r$ is, as $k'_h > k_h$.
Since $B^{k_r}_r= B^{(k_r)}_{(r)}$ the claim is  proved.

By the claim, we necessarily have
\begin{equation}
\label{eq-r*}
 B_{(*)}  \cap C^{}_* = \emptyset. 
 \end{equation} 
 If $h<r$   a stronger fact holds, to state which 
  we depart momentarily from the notation in Proposition~\ref{hypres} and  write $S_{k_h}\to  S^1_{h+1}$ for the blow-up at $q^1_{h+1}$; denote by $C^{k_h}, B^{(k_h)}\subset S_{k_h}$ the proper transforms of $B$ and $C$.
By the claim,  $C^1_{h+1}$ is not tangent   to $B^{(1)}_{(h+1)}$, hence 
\begin{equation}
\label{eq-kh}
C^{k_h} \cap B^{(k_h)}= \emptyset  
 \end{equation} 
 (which is stronger than \eqref{eq-r*}).

We are ready to get a contradiction. Assume  $h< r$;  by iterating  \eqref{eq-int}  we have
\begin{equation}
\label{eq-doth}  (C^{k_h} \cdot  B^{(k_h)}) =  bd-\sum_{j=1}^{h}k_jl_jl'_j + l_{h+1}l'_h. 
 \end{equation} 
Indeed,    $C^i_j$ has an $l_j$-fold point at $q^i_j$ for   $j=1,\ldots, h$ and $i=1,\ldots, k_h$,
  while $C^1_{h+1}$ has an $l_{h+1}$-fold point. The curve  $B^{(i)}_{(j)}$ has a $l'_j$-fold point at $q^i_j$ for   $j=1,\ldots, h$ and $i=1,\ldots, k_j$, and $B^{(1)}_{(h+1)}$ has an $l'_h$-fold point. Hence  \eqref{eq-doth} follows.
We  combine  \eqref{eq-doth} with  Lemma~\ref{lm-key} (which can be applied as $k_j=k'_j$ for all $j<h$) and  we get 
 $$(C^{k_h} \cdot  B^{(k_h)})>bd-nn'\geq 0 $$ 
 (as $b\geq n'$ and $d\geq n$),
 contradicting \eqref{eq-kh}.
  
Assume   $h=r$, hence $k_j=k_j'$ for all $j<r$. Then for   $j=1,\ldots, r$ and $i=1,\ldots, k_r$, at $q^i_j$
     the curve  $C^i_j$ has an $l_j$-fold point, and 
 $B^{(i)}_{(j)}$   has 
$l'_j$-fold point.  Therefore, 
\begin{equation}
\label{eq-C*}
   (C^{}_* \cdot B_{(*)}) = bd-\sum_{j=1}^{r}k_jl_jl'_j.  
\end{equation}
Since $l_{r+1}=0$ we get
$$
\sum_{j=1}^{r}k_jl_jl'_j=\sum_{j=1}^{r}k_jl_jl'_j +l_{r+1}l'_{r}<nn'
$$
by Lemma~\ref{lm-key} with $h=r$. 
  Hence    $$
  (C^{}_* \cdot B_{(*)})  >bd-nn'\geq 0. 
$$ 
contradicting \eqref{eq-r*}.    The theorem is proved.
\end{proof}

We used the following elementary Lemma. 

\begin{lemma} 
\label{lm-key}  
Let   \eqref{eq-eu} and \eqref{eq-eu'} be the     euclidean 
sequences of $(m,n)$ and $(m',n')$. Assume  $ r\leq r'$ and let $h\leq r$. If  $k_j=k'_j$ for every $j=1,\ldots h-1$ then 
$$
\sum_{j=1}^{h}k_jl_jl'_j +l_{h+1}l'_{h}<nn'.
$$

\end{lemma}
\begin{proof}
The proof is based on the following   identity:
 \begin{equation}
\label{eq-key}
\sum_{j=1}^{h}k_jl_jl'_j +l_{h+1}l'_{h}= \left\{ \begin{array}{ll}
l_1l_0' &  \  \text {if $h$ is even}  \\
l_0l_1'  &  \    \text {if $h$ is odd.} 
 \end{array} \right.
\end{equation}
We have $m=l_1<n=l_0$ and  $m'=l'_1<n=l'_0$,
hence $ l_1l_0' <nn'$   and $ l_1'l_0 <n'n$.
Combining this with  \eqref{eq-key}  the Lemma is proved.

We prove   \eqref{eq-key}  by induction on $h$,    applying a few times the  basic identities:
$$
  l_{j-2} =k_{j-1}l_{j-1}+ l_j,\quad \quad  l'_{j-2} =k'_{j-1}l'_{j-1}+ l'_j.
$$
The base cases are $h=1$ and $h=2$. If $h=1$   we get 
   $$
   k_1l_1l'_1+l_2l'_1=(k_1 l_1+ l_2)l'_1 =l_0l_1'.
   $$
   If $h=2$ then $k_1=k'_1$, and we get
    $$
   k_1l_1l'_1+k_2l_2l'_2+l_3l'_2= k'_1l_1l'_1+(k_2l_2+l_3)l'_2= k'_1l_1l'_1+l_1l'_2=l_1(k'_1l'_1+ l'_2)=l_1l_0'. 
   $$
   Let  $h\geq 3$. If  $h$ is even  we have
\begin{align*}
&\sum_{j=1}^{h}k_jl_jl'_j+l_{h+1}l'_{h} =\sum_{j=1}^{h-1}k_jl_jl'_j+ k_{h}l_{h}l'_{h}+l_{h+1}l'_{h}= \\
=&\sum_{j=1}^{h-1}k_jl_jl'_j+ ( k_{h}l_{h} +l_{h+1})l'_{h} = \sum_{j=1}^{h-1}k_jl_jl'_j+ l_{h-1}l'_{h} = \\
=&\sum_{j=1}^{h-2}k_jl_jl'_j+l_{h-1}( k_{h-1} l'_{h-1}+ l'_{h}) =\sum_{j=1}^{h-2}k_jl_jl'_j+l_{h-1}( k'_{h-1} l'_{h-1}+ l'_{h}) = \\
=&\sum_{j=1}^{h-2}k_jl_jl'_j+l_{h-1}l'_{h-2} =l_1l_0' \quad \quad  \quad \quad  \quad \quad 
\end{align*}
  by the induction hypothesis ($h-2$ is even). 
 The case $h$ odd follows in the same way. The Lemma is proved.
  \end{proof} 
  
  The following special case   of the theorem is worth an explicit mention.
  \begin{cor}
\label{cor-mirror}
Hypotheses of  Theorem~\ref{mirror}.
 
If  $\gcd(m,n)=\gcd(m',n')$ then $(m,n)=(m',n')$; in particular, if $B$ and $C$ are smooth at $q$, then $n=n'$.
 \end{cor}

\begin{example}
 Let $B  \subset \PP^2$ be the projective conic of affine equation $y=x^2$,  and $C\subset \PP^2$ be the projective    cubic of affine equation $y=x^2+y^3$. A simple computation shows that $B$ and $C$ intersect only at $q=(0,0)$. The theorem implies that $q$ is not a flex of $C$  (which  is easy to check in this case).
\end{example}

  \subsection{Intersection multiplicity at equivalent points}
  
We proved in
  Theorem~\ref{mirror}  that if two curves $B$ and $C$ are hypertangent in one point $q$, then $q\in B$ is equivalent to $q\in C$, i.e. their contact tropical curves are isomorphic. It is  natural to pose the opposite problem, namely, what can be said about the intersection of   two curves   at equivalent points.

  \begin{prop}
 \label{inverse}
 Let $B,C\subset \PP^2$ be integral curves  and $q\in B\cap C$.
Assume that $q$ is  an $(m',n')$-point of $B$ and an $(m,n)$-point of $C$, with  $(m,n)\sim(m',n')$; set $k=\lfloor n/m\rfloor$.
If  $(B\cdot C)_q>m'm$ (i.e. $B$ and $C$ are tangent in $q$), then
\begin{equation}
  (B\cdot C)_q\geq \left\{ \begin{array}{ll}
km'm &  \  \text {if  }  n/m\in \N \\
km'm +1+(n'-m')(n-m)&  \  \text {if  } n/m\not \in \N . 
 \end{array} \right.
\end{equation}
\end{prop}
\begin{proof} Let $L$ be the tangent line at $q$ to  $C$ and $B$, and 
 $S_1\to \PP^2$  the blow-up at $q$; denote by $L^1,B^1,C^1\subset S_1$ the proper transforms of  $L,B,C$, and by $q^1\in L^1$ the point lying over $q$. 
Denote by $S_2\to S_1$ the blow-up at $q^1$, by $L^2,B^2,C^2\subset S_2$ the proper transforms of  $L,B,C$, and by $q^2\in L^2$ the point lying over $q^1$. Iterating we get a sequence of blow-ups
 $$S_{k+1}\to S_k\to S_{k-1}\to\ldots  S_{i+1}\to S_i\to \ldots \to S_1\to S_0=\PP^2$$ such that $L^i,B^i,C^i\subset S_i$ are the proper transforms of 
 $L,B,C$ and $S_{i+1}\to S_i$ is the blow-up at $q^i\in L^i$, with $q^i$ lying over $q$; we write $E_i\subset S_i$ for the exceptional divisor.
 
By Lemma~\ref{lm-k}, for every $i\leq k-1$ the curve  $C^i$, respectively   $B^i$, has a point of multiplicity $m$, respectively   $m'$, lying over $q$. Therefore, setting
 $$(C^i\cdot L^i)_{E_i}=\sum_{p\in E_i}(C^i\cdot L^i)_p$$
  we have, as $k\leq n/m$
 $$
 (C^i\cdot L^i)_{E_i}=(C\cdot L)_q-im=n-im\geq n-(k-1)m\geq n-(n/m-1)m=m\geq 1
 $$
 hence, as $L^i\cap E_i= q^i$, we get
$(C^i\cdot L^i)_{q^i}\geq 1$. An analogous argument applied to $B$ gives 
$(B^i\cdot L^i)_{q^i}\geq 1$ for all $i\leq k-1$. Therefore  $ q^{k-1}\in B^{k-1}\cap C^{k-1}$ hence 
$$
(B^{k-1}\cdot C^{k-1})_{q^{k-1}}\geq m'm.
$$
On the other hand 
$$
(B^{k-1}\cdot C^{k-1})_{q^{k-1}}=(B\cdot C)_q-(k-1)m'm
$$
hence
$$
(B\cdot C)_q\geq m'm+(k-1)m'm=km'm.
$$
If $n/m\in \N$ there is nothing left to prove.
Let   $n/m\not\in \N$, hence  $k<n/m$. 
Then  $C^k$  has an 
$(n-m,m)$-point w.r.t. $E_k$, and   $B^k$ has an $(n'-m',m')$-point w.r.t. $E_k$   (by Lemma~\ref{lm-k}); let us show that this point is $q^k$ for both $C^k$ and $B^k$.
We have
$$
(C^k\cdot L^k)_{q^{k}}=n-km>n-nm/m=0
$$
hence $q^k\in C^k$. The same argument shows that $q^k\in B^k$. Hence $C^k$ and $B^k$ intersect in $q^k$,
as wanted.  Now, since   $C^k$ and $B^k$ are both tangent to $E_k$, their intersection is not transverse, hence their proper transforms in $S_{k+1}$ satisfy
$$(B^{k+1}\cdot C^{k+1})_{{\tilde q}^{k+1}}\geq 1$$
with ${\tilde q}^{k+1}$ lying over $q$.  On the other hand
$$(B^{k+1}\cdot C^{k+1})_{{\tilde q}^{k+1}}=(B\cdot C)_q-km'm-(n'-m')(n-m)
$$
hence
$ 
(B\cdot C)_q\geq 1+km'm+(n'-m')(n-m)
$. The proof is complete.
\end{proof}

\section{Genus and dimension formulas}
In this section we obtain bounds for the delta invariant (defined below)  of a unibranch point of a   curve $C \subset \PP^2$, and  for the codimension of the locus of   plane curves with a $(m,n)$ point.
 
\subsection{Genus computation}
 Let $C\subset \PP^2$ and $q\in C$; let $\deg C=d\geq 2$. 
The  geometric genus, $g(C)=g(C^{\nu})$ (with $C^{\nu}$   the normalization of $C$),   satisfies  
$$
g(C)\geq (d-1)(d-2)/2-\delta (q)
$$
with equality if and only if $q$ is the only singular point of $C$. The number $\delta  (q)$ is called the {\it delta-invariant} of $q$.
 The following formula is well-known
\begin{equation}
 \label{eq-di}
  \delta(q)=\sum m_p(m_p-1)/2
\end{equation}
 where $p$ varies over all points infinitely near to $q$ (the points added over $q$ in the proper transforms of the blow-ups), and $m_p$ is the multiplicity of $p$. For example, if $q$ is an ordinary $m$-fold point (i.e. $C$ has   $m$ branches meeting transversally in $q$), then $  \delta(q)=m(m-1)/2$.

\begin{prop}
\label{delta}  
 Let  $C\subset \PP^2$ have degree $d\geq 2$ and  let  $q\in C$ be  an  $(m,n)$-point; set $\gcd (n,m)=c$.  
Then   $$
\delta  (q)\geq  (nm - n - m + c)/2 
$$
with equality if   $c=1$.
\end{prop}
\begin{proof}
 Proposition~\ref{hypres} considers
the  desingularization, $\nu_q:C^{\nu}_q\to C$, of $C$ at $q$, and describes
  the points infinitely near to $q$ appearing in the partial desingularization given by $\beta$.
  This, combined with 
 \eqref{eq-di},  
gives  
\begin{equation}
 \label{eq-dq}
 \delta (q)\geq  \sum_{j=1}^{r}k_j l_j(l_j-1)/2.  
\end{equation}
 Indeed, for every $j=1,\ldots, r$  we have exactly $k_j$ infinitely near points of multiplicity $l_j$, 
 and each such point contributes by $l_j(l_j-1)/2$. 
 Note that if $c=1$   the summand corresponding to $j=r$ vanishes, as $l_{r}=1$. In this case
   the desingularization of $C$ at $q$ is actually included in $\beta$,
 hence there are no other infinitely near points, therefore in \eqref{eq-dq} we have equality.

 The following Lemma~\ref{contoLaura} completes the proof.
\end{proof}

 \begin{lemma}\label{contoLaura} Let $n>m>0$ be two integers and consider the euclidean sequence \eqref{eq-eu}. Then   $k_j = (l_{j-1} - l_{j+1})/l_j$ and
  \begin{equation}\label{eq-Laura}
   \sum_{j=1}^{r} k_j l_j (l_j-1)/2 = (nm -n - m + c)/2.  
  \end{equation}
 \end{lemma}
\begin{proof}
The expression of $k_j$ follows directly from its definition and the construction of the euclidean sequence.
The second assertion follows from a special case of the following formula
\begin{equation}\label{eq-Laura2}
 \sum_{j=1}^{r} (l_{j-1} - l_{j+1})(l_j - 1) = l_0l_1 - l_0 -l_1 + l_{r+1} + l_{r} - l_{r}l_{r+1}. 
\end{equation}
In our case $l_0 = n$, \  $l_1 = m$,\  $l_{r} = c$  and $l_{r+1}=0$,   from which \eqref{eq-Laura} follows. 
To prove   \eqref{eq-Laura2} we proceed by induction on $r$. If $r=1$ then
\[
 (l_0 - l_2)(l_1 - 1) = l_0l_1 - l_1l_2 - l_0 + l_2,
\]
as wanted. To prove the general case, the induction hypothesis reads as follows
\begin{equation*}
 \sum_{j=1}^{r-1} (l_{j-1} - l_{j+1})(l_j - 1) = l_0l_1 - l_0 -l_1 + l_{r-1} + l_{r} - l_{r}l_{r-1}. 
\end{equation*}
Thus we obtain
\begin{align*}
 \sum_{j=1}^{r} (l_{j-1} - l_{j+1})(l_j - 1)  = \\
 l_0l_1 - l_0 -l_1 + l_{r-1} + l_{r} - l_{r}l_{r-1} +(l_{r-1} - l_{r+1})(l_{r} - 1) = \\
 l_0l_1 - l_0 -l_1 + l_{r-1} + l_{r} - l_{r}l_{r-1} +l_{r-1}l_{r} - l_{r}l_{r+1} - l_{r-1} + l_{r+1}    = \\
 l_0l_1 - l_0 -l_1 + l_{r}- l_{r}l_{r+1 }+ l_{r+1}.
\end{align*}
The Lemma is proved.    \end{proof}

\subsection{Dimension computation}
Let $n,m$ be two integers with $n>m>0$; 
let $L\subset \PP^2$ be a line and $q\in L$ a point.
We denote by  
 $$U^{(m,n)}_d(L;q)\subset  \PP^{d(d+3)/2}$$ 
 or simply by $U^{(m,n)}_d$, 
  the space of integral curves of degree $d\geq 2$ having a unibranch  $(m,n)$-point at $q$ with tangent line equal to $L$. As it will be clear   (if it is not already),   its closure,
 $\ov{U^{(m,n)}_d }$, is a linear subspace of $\PP^{d(d+3)/2}$.
  We write  $ \codim U^{(m,n)}_d =\codim \ov{U^{(m,n)}_d }$  for the  codimension in $\PP^{d(d+3)/2}.$

\begin{prop}
\label{hyptan} Assume $d\geq n>m>0$ and  set $\gcd (n,m)=c$.  
Then  
  $$
 \codim U^{(m,n)}_d \geq (nm+m+n+c-2)/2 
$$
with equality when $c=1$.
\end{prop}
 
\begin{proof} We choose, as usual, affine coordinates  $(x,y)$ so that
  $L$ has equation $y=0$ and $q=(0,0)$. Let
  $C\in  U^{(m,n)}_d(L;q)$ 
and let its  affine equation be $f(x,y)=0$ with
$$ 
f(x,y)=\sum_{0\leq i+j\leq d}a_{i,j}x^iy^j.
$$ 
By hypothesis  $a_{n,0}, a_{0,m}\neq 0$.

As $(C\cdot L)_q=n$ 
we have   $f(x,0)= x^np(x)$  with $p(0)\neq 0$, hence the smallest power of $x$ appearing in $f$ is $x^n$. This gives   the following  $n$  conditions
\begin{equation}
 \label{eq-1}
 a_{i,0}=0 \quad \forall i=0,\ldots, n-1.
\end{equation}
 Next,  as $q $  is an $m$-fold point we have
 $ 
 a_{i,j}=0$ for all $i+j\leq m-1$;  these are 
  $m(m+1)/2$ conditions, but exactly $m$ of them 
  are also in \eqref{eq-1}. Hence  the number of  new independent conditions on  $f$ is equal to
 \begin{equation}\label{eq-2} 
m(m+1)/2-m=m(m-1)/2.
\end{equation} 

Next, $C$ is unibranch at $q$, and its tangent line there is  $y=0$. Therefore the only summand of degree $m$   is $a_{0,m}y^m$ and  we get the following  new  $m-1$ conditions:
 \begin{equation}\label{eq-3}
 a_{i,j}=0,\quad \forall i+j=m, \quad i,j\geq 1.
\end{equation}

Now  consider the blow-up of $\PP^2$ at $q$,   let  $C'$ and $L'$ be the strict transforms of $C$ and $L$, and $q'\in C'$ the point lying over $q$.
We set $y=vx$ and use $(x,v)$ as affine coordinates in the blow-up.   The   equation of $C'$ is   
$$
f'(x,v)=a_{0,m}v^m+\sum_{m+1\leq i+j\leq d} a_{i,j}x^{i+j-m}v^j=a_{0,m}v^m+\sum_{\stackrel{1\leq h\leq d-m}{m+1\leq l\leq d}} a_{h-l+m,l}x^hv^l 
$$
with $j=l$ and $i=h-l+m$.
The   conditions imposed earlier imply exactly that  $f'$ has no summand of type $x^h$ for all $h=0,\ldots, d-m-1$, and no summand of type $v^l$ for all $l=1,\ldots, m-1$.  

Let $m'$ be the multiplicity of $q'\in C'$;  the fact that $C'$ is unibranch   at $q'$   yields  new conditions; 
we claim that the  number  of  them     is  
$$m'(m'-1)/2 .$$

If $n<2m$ we have   $n-m=m'$ and $q'$ is a $(m',m)$-point of $C'$ with respect to exceptional divisor, by Lemma~\ref{lm-k}. More precisely,
the  term of smallest degree  of $f'$ is    $a_{n,0}x^{m'}$, and  $C'$ has a unibranch  $m'$-fold point  at $q'$.  Hence every   coefficient  of a term  of degree less than ${m'}$, and of 
  a term  of degree equal to ${m'}$  other than $x^{m'}$, must vanish.
These are ${m'}({m'}+1)/2+{m'}$ independent conditions, but  some of them are already satisfied. Indeed, we   know that 
 $f'$ has no term of type $x^l$ for all $l=0,\ldots, {m'}-1$ ($m'$ conditions), and no term of type $v^n$ for all $n=1,\ldots, {m'}$ (${m'}$ conditions).  There are no other cases, as if $k,l>0$ the coefficient of  $x^kv^l$ is $a_{k-l+m,l}$ corresponding to a summand in $f$ of degree $k+m\geq m+1$, and no conditions had been imposed on such coefficients.  
Hence  the number of new conditions   is  
\begin{equation}
 \label{eq-r}
{m'}({m'}+1)/2+{m'}-{m'}-{m'}={m'}({m'}-1)/2. 
\end{equation}
Let $n> 2m$.    Then $m=m'$ and    $C'$ has a an $(m,n-m)$-point at $q'$ with respect to $L'$, by Lemma~\ref{lm-k}.
Hence in $f'$ there are no monomials
 of degree less than $m$, and  no  monomials 
  of degree   $m$,  except   $v^m$.
  These are $m(m+1)/2+m$   conditions,   some of which  are already satisfied, indeed
 $f'$ has no term of type $x^l$ for all $l=0,\ldots, m$ ($m+1$ conditions), and no summand of type $v^n$ for all $n=1,\ldots, m-1$ ($m-1$ conditions).  
Hence  the number of new conditions   is  
\begin{equation}
 \label{eq-4}
m(m+1)/2+m -(m+1)-(m-1) =m(m-1)/2. 
\end{equation}
Let $n= 2m$.    
Then  $C'$ has a an $m$-fold point at $q'$.
Hence in $f'$ there are no monomials
 of degree less than $m$. Moreover,  as we saw in the proof of Lemma~\ref{lm-k}, the homogeneous part of degree $m$ of $f'$ has the form $(\alpha x+\beta v)^m$ with $\alpha^m=a_{2m,0}$ and $\beta^m=a_{0,m}$. This imposes $m-1$ conditions. 
 We get a total of $m(m+1)/2+m-1$ conditions, but, as in the earlier cases, some of them already hold.
In fact
 $f'$ has no term of type $x^l$ for all $l=0,\ldots, m-1$ ($m$ conditions), and no summand of type $v^n$ for all $n=1,\ldots, m-1$ ($m-1$ conditions).  
Hence  the number of new conditions   is  
\begin{equation}
 \label{eq-5}
m(m+1)/2+m-1 -m-(m-1)=m(m-1)/2.  
\end{equation}
The claim is proved.

Now we need to apply the claim to   all  the curves appearing in the chain    described in Proposition~\ref{hypres}, of which we use the notation. Proposition~\ref{hypres} gives a sequence of curves whose multiplicity at the point lying over $q$ is equal to    $l_j$, with $j=1,\ldots,r$. By the claim, each of these  gives  $l_j(l_j-1)/2$  new conditions.
If $j\geq 2$ we get  exactly $k_j$ such curves, hence a total of 
  $k_j l_j(l_j-1)/2$ conditions.
 For $j=1$ we   have only $k_1-1$ such  blow-ups, hence  $(k_1-1)m(m-1)/2$ conditions. But in the first part of the proof we had   $m(m-1)/2$ conditions; see \eqref{eq-2}.  Hence also for $j=1$ we have
 $k_1 m(m-1)/2$ conditions, hence
 $$
\sum_{j=1}^{r}k_j l_j(l_j-1)/2
$$ 
conditions. 
Adding these to the ones computed in \eqref{eq-1} and \eqref{eq-3} we
get
$$ \codim U^{(m,n)}_d= 
 n+(m-1) +\sum_{j=1}^{r}k_j l_j(l_j-1)/2.$$ 
By applying 
  Lemma~\ref{contoLaura} we are done.
   \end{proof}

\bibliography{references}{}
\bibliographystyle{alpha}
\end{document}